\documentstyle[12pt]{article}
\title{Spectral Conditions for the Bipancyclic Bipartite Graphs}
\author {Rao Li \\
         Dept. of mathematical sciences \\
         University of South Carolina Aiken \\
	   Aiken, SC 29801 \\
	   USA \\
         {\it Email: raol@usca.edu }
         }
\date{July 28, 2020}

\begin{document}
\maketitle
\begin{abstract}
Let $G = (X, Y; E)$ be a bipartite graph with two vertex partition subsets  $X$ and $Y$. 
 $G$ is said to be balanced if $|X| = |Y|$.
 $G$ is said to be bipancyclic if it contains cycles of every even length from $4$ to $|V(G)|$.
 In this note, we present spectral conditions for the bipancyclic bipartite graphs. 
 \end{abstract} 
$$Keywords: Spectral \,\, Condition, \,\,\, Bipancyclic \,\, Bipartite \,\, Graphs.$$ 
$$2010 \, Mathematics \, Subject \, Classification: 05C50; \, 05C45 $$ \\

\noindent {\bf 1.  Introduction} \\

We consider only finite undirected graphs without loops or multiple edges.
Notation and terminology not defined here follow that in \cite{Bondy}.
Let $G = (V(G), \, E(G))$ be a graph. The graph $G$ is said to be Hamiltonian if it contains a cycle of length $|V(G)|$.
The graph $G$ is said to be pancyclic if it contains cycles of every length from $3$ to $|V(G)|$.
Let $G = (V_1(G), V_2(G); \, E(G))$ be a bipartite graph with two vertex partition subsets  $V_1(G)$ and $V_2(G)$. 
The bipartite graph $G$ is said to be semiregular bipartite if all the vertices in $V_1(G)$ have the same degree and 
all the vertices in $V_2(G)$ have the same degree.
 The bipartite graph $G$ is said to be balanced if $|V_1| = |V_2|$. Clearly, if a bipartite graph is Hamiltonian, then it must be balanced.
 The bipartite graph $G$ is said to be bipancyclic if it contains cycles of every even length from $4$ to $|V(G)|$.
 The balanced bipartite graph $G_1 = (A, B; E)$ of order $2n$ with $n \geq 4$ is defined as follows:
 $A = \{a_1, a_2, ... , a_n\}$, $B = \{b_1, b_2, ... , b_n\}$, and 
 $E = \{a_ib_j : 1 \leq i \leq 2, (n - 1) \leq j \leq n \} \cup \{a_ib_j : 3 \leq i \leq n, 1 \leq j \leq n \}$. 
 Notice that $G_1$ is not Hamiltonian.  \\
 
 The eigenvalues of a graph $G$, denoted $\lambda_1(G) \geq \lambda_2(G) \geq \cdots \geq \lambda_n(G)$,
 are defined as the eigenvalues of its adjacency matrix $A(G)$.
 Let $D(G)$ be a diagonal matrix such that its diagonal entries are the degrees of vertices in a graph $G$.
 The Laplacian matrix of a graph $G$, denoted $L(G)$, is defined as
$D(G) - A(G)$, where $A(G)$ is the adjacency matrix of $G$. The eigenvalues 
$\mu_1(G) \geq \mu_2(G) \geq \cdots \geq \mu_{n - 1}(G) \geq \mu_{n}(G) = 0$ of $L(G)$ are called
the Laplacian eigenvalues of $G$. The second smallest Laplacian eigenvalue $\mu_{n - 1}(G)$ 
is also called the algebraic connectivity of the graph $G$ (see \cite{F}). 
The signless Laplacian matrix of a graph $G$, denoted $Q(G)$, is defined as
$D(G) + A(G)$, where $A(G)$ is the adjacency matrix of $G$. The eigenvalues 
$q_1(G) \geq q_2(G) \geq \cdots \geq q_n(G)$ of $Q(G)$ are called
the signless Laplacian eigenvalues of $G$.\\
 
Yu et al. in \cite{Yu} obtained some spectral conditions for the pancyclic graphs.  Motivated by the results in \cite{Yu}, 
we present spectral conditions for the bipancyclic bipartite graphs. The main results are as follows. \\

\noindent {\bf Theorem $1$}. Let $G = (X, Y; E)$ be a connected balanced bipartite graph of order $2n$ with $n \geq 4$, $e$ edges, and $\delta \geq 2$.
If $\lambda_1 \geq \sqrt{n^2 - 2n + 4}$, then $G$ is bipancyclic. \\

\noindent {\bf Theorem $2$}. Let $G = (X, Y; E)$ be a connected balanced bipartite graph of order $2n$ with $n \geq 4$, $e$ edges, and $\delta \geq 2$. If
 $$\mu_{n - 1} \geq \frac{2(n^2 - 2n + 4)}{n},$$ 
 then $G$ is bipancyclic. \\

\noindent {\bf Theorem $3$}. Let $G = (X, Y; E)$ be a connected 
balanced bipartite graph of order $2n$ with $n \geq 4$, $e$ edges, and $\delta \geq 2$. If 
$$q_1 \geq \frac{2(n^2 - n + 2)}{n},$$ then $G$ is bipancyclic.\\ 

\noindent {\bf 2.  Lemmas} \\

In order to prove the theorems above, we need the following results as our lemmas. \\

Lemma $1$ is the main result in \cite{MS}, \\

\noindent {\bf Lemma $1$}. Let $G = (X, Y; E)$ be a balanced bipartite graph of order $2n$ with $n \geq 4$. Suppose that
 $X = \{x_1, x_2, ... , x_n\}$, $Y = \{y_1, y_2, ... , y_n\}$, 
$d(x_1) \leq d(x_2) \leq \cdots \leq d(x_n)$, and  
$d(y_1) \leq d(y_2) \leq \cdots \leq d(y_n)$. 
If $$d(x_k) \leq k < n \Longrightarrow d(y_{n - k}) \geq n - k + 1,$$
then $G$ is bipancyclic. \\

Lemma $2$ below follows from Proposition $2.1$ in \cite{Fri}. \\ 

\noindent {\bf Lemma 2.} Let $G$ be a connected bipartite graph of order $n \geq 2$ and $e \geq 1$ edges. Then $\lambda_1 \leq \sqrt{e}$. 
If $\lambda_1 = \sqrt{e}$,  then $G$ is a complete bipartite graph $K_{s, \, t}$, where $e = s \, t$. \\

Lemma $3$ below is $4.1$ in \cite{F}. \\

\noindent {\bf Lemma 3.} Let $G$ be a noncomplete graph. Then $\mu_{n - 1} \leq \kappa$, where $\kappa$ is the vertex connectivity of $G$. \\

Lemma $4$ below is Theorem $2.9$ in \cite{Bo}. \\

\noindent {\bf Lemma 4.} Let G be a balanced bipartite graph of order $2n$ and $e$ edges. Then
$q(G) \leq \frac{e}{n} + n$. \\

Lemma $5$ below is Lemma $2.3$ in \cite{Feng}. \\

\noindent {\bf Lemma 5.} Let $G$ be a connected graph. Then
$$q_1 \leq \max\{ d(u) + \frac{\sum_{v \in N(u)} d(v)} {d(u)} : u \in V \},$$ 
with equality holding if and only if $G$ is either semiregular bipartite or regular. \\

\noindent {\bf Lemma 6.} Let $G = (X, Y; E)$ be a balanced bipartite graph of order $2n$ with $n \geq 4$, $e$ edges, and $\delta \geq 2$.
If $e \geq n^2 - 2n + 4$, then $G$ is bipancyclic or $G = G_1$. \\

\noindent {\bf Proof of Lemma 6.} Without loss of generality, we assume that  $X = \{x_1, x_2, ... , x_n\}$, $Y = \{y_1, y_2, ... , y_n\}$, 
$d(x_1) \leq d(x_2) \leq \cdots \leq d(x_n)$, 
and $d(y_1) \leq d(y_2) \leq \cdots \leq d(y_n)$. Suppose $G$ is not bipancyclic.  
Then Lemma $1$ implies that there exists an integer $k$ such that $1 \leq k < n$, $d(x_k) \leq k$, and $d(y_{n - k}) \leq n - k$.  Thus
$$2n^2 - 4n + 8 \leq 2 e = \sum_{i = 1}^n d(x_i) + \sum_{i = 1}^n d(y_i) $$
$$\leq k^2 + (n - k) n + (n - k)^2 + k n $$
$$= 2n^2 - 4n + 8 - (k - 2)(2n - 2k - 4).$$
Since $\delta \geq 2$, we have that $k \neq 1$. Therefore we have the following possible cases. \\

\noindent{\bf Case $1$.} $k = 2$. \\

In this case, all the inequalities in the above arguments now become equalities. Thus $d(x_1) = d(x_2) = 2$,
$d(x_3) = \cdots = d(x_n) = n$, $d(y_1) = \cdots = d(y_{n - 2}) = n - 2$, and $d(y_{n - 1}) = d(y_{n - 2}) = n$.
Hence $G = G_1$. \\

\noindent{\bf Case $2$.} $(2n - 2k - 4) = 0$. \\

In this case, we have $n = k + 2$ and all the inequalities in the above arguments now become equalities. Thus $d(x_1) = \cdots = d(x_{n - 2}) = n - 2$,
$d(x_{n - 1}) = d(x_n) = n$, $d(y_1) = d(y_{2}) = n - 2$, and $d(y_{3}) = \cdots = d(y_{n - 2}) = n$.
Hence $G = G_1$. \\

\noindent{\bf Case $3$.} $k \geq 3$ and $2n - 2k - 4 < 0$. \\

In this case, we have that $n < k + 2$, namely, $n \leq k + 1$. Since $k < n$, we have $k = n - 1$. This implies that $d(y_1) \leq 1$, contradicting to 
the assumption of $\delta \geq 2$.  \\

This completes the proof of Lemma $6$. \hfill{$\Box$} \\

\noindent {\bf 3.  Proofs} \\

\noindent {\bf Proof of Theorem $1$}. Let $G$ be a graph satisfying the conditions in Theorem $1$. Then we, from Lemma $2$, have that 
$$\sqrt{n^2 - 2n + 4} \leq \lambda_1 \leq \sqrt{e}.$$ 
Thus $e \geq  n^2 - 2n + 4$. Therefore by Lemma $6$ we have that $G$ is bipancyclic or $G = G_1$. \\

If  $G = G_1$, then $e = n^2 - 2n + 4$.
Hence  $$\sqrt{n^2 - 2n + 4} \leq \lambda_1 \leq \sqrt{e} = \sqrt{n^2 - 2n + 4}.$$
So $\lambda_1 = \sqrt{e}$. Lemma $2$ implies that $G_1$ is a complete bipartite graph, a contradiction. \\

This completes the proof of Theorem $1$. \hfill{$\Box$} \\

\noindent {\bf Proof of Theorem $2$}. Let $G$ be a graph satisfying the conditions in Theorem $2$. Then we, from Lemma $3$, have that 
$$\frac{2(n^2 - 2n + 4)}{n} \leq \mu_{n - 1} \leq \kappa \leq \delta \leq \frac{2e}{n}.$$ 
Thus $e \geq  n^2 - 2n + 4$. Therefore by Lemma $6$ we have that $G$ is bipancyclic or $G = G_1$. \\

If  $G = G_1$, then $e = n^2 - 2n + 4$. 
Hence 
$$\frac{2(n^2 - 2n + 4)}{n} \leq \mu_{n - 1} \leq \kappa \leq \delta \leq \frac{2e}{n} = \frac{2(n^2 - 2n + 4)}{n}.$$ 
This implies that $G_1$ is a regular graph, a contradiction. \\

This completes the proof of Theorem $2$. \hfill{$\Box$} \\

\noindent {\bf Proof of Theorem $3$}. Let $G$ be a graph satisfying the conditions in Theorem $3$. Then we, from Lemma  $4$, have that 
$$\frac{2(n^2 - n + 2)}{n} \leq q_1 \leq \frac{e}{n} + n.$$
Thus $e \geq  n^2 - 2n + 4$. Therefore by Lemma $6$ we have that $G$ is bipancyclic or $G = G_1$. \\

If $G = G_1$, then $e = n^2 - 2n + 4$. Therefore
$$\frac{2(n^2 - n + 2)}{n} \leq q_1 \leq \frac{e}{n} + n = \frac{ n^2 - 2n + 4}{n} + n = \frac{2(n^2 - n + 2)}{n}.$$
Hence $$q_1 = \frac{2(n^2 - n + 2)}{n}.$$ 
It can be verified that 
$$\max\{ d(u) + \frac{\sum_{v \in N(u)} d(v)} {d(u)} : u \in V \} = \frac{2(n^2 - n + 2)}{n} = q_1.$$ 
Thus Lemma $5$ implies that $G_1$ is semiregular or regular, a contradiction. \\ 

 This completes the proof of Theorem $3$. \hfill{$\Box$} \\


\begin{thebibliography}{20}
  

 \bibitem{Fri} A. Bhattacharya, S. Friedland, and U. Peled, On the first eigenvalue of bipartite graphs, 
 The Electronic Journal of Combinatorics 15 (2008), \#R 144.  
 
 \bibitem{Bondy} J. A. Bondy and U. S. R. Murty, Graph Theory with
          Applications, Macmillan, London and Elsevier, New York (1976).  
          
\bibitem{C} D. Cvetkovi\'{c}, M. Doob, and H. Sachs, Spectra of Graphs, 3rd Edition,
         Johann Ambrosius Barth (1995).
         
\bibitem{Feng} L. Feng and G. Yu, On three conjectures involving the signless Laplacian spectral
radius of graphs, Publ. Inst. Math. 85 (2009), 35--38.
                  
 \bibitem{F} M. Fiedler, Algebraic connectivity of graphs,   Czechoslovak Mathematical Journal 23 (1973), 298--305.        
          
 \bibitem{Bo} B. Li and B. Ning, Spectral analogues of Erd\H{o}s' and Moon-Moser's theorems on
                        Hamilton cycles, https://arxiv.org/pdf/1504.03556.pdf, retrieved on July 28, 2020.
                                                                  
\bibitem{MS}  J. Mitchem and E. F. Schmeichel, Bipartite graphs with all even lengths, Journal of Graph Theory 6 (1982), 429--439.   
          
\bibitem{Yu}  G. Yu, T. Yu, A. Shu, and X. Xia, Some sufficient conditions on pancyclic graphs, 
                            https://arxiv.org/pdf/1809.09897.pdf, retrieved on July 28, 2020.  
                            



\end{thebibliography}
\end{document}